\documentclass[12pt]{article}
\usepackage{amsmath,amssymb,amsthm}
\usepackage{graphics,epsfig}
\usepackage{hyperref}
\usepackage{natbib}
\usepackage{color}
\usepackage{graphicx}
\usepackage{caption}
\usepackage{subcaption}
\usepackage{float}
\usepackage{mathrsfs}
\usepackage{multirow}
\usepackage{comment}
\usepackage{setspace}
\usepackage{bbm}
\usepackage{bm}
\usepackage{amsfonts}
\usepackage{xcolor}
\usepackage{pslatex}
\usepackage{setspace}

\begin{document}

  \title{\bfseries Path integral control under McKean-Vlasov dynamics}

  \author{
Timothy Bennett\footnote{Corresponding author, {\small\texttt{tdb1923@jagmail.southalabama.edu}}}\; \footnote{Department of Mathematics and Statistics,  University of South Alabama, Mobile, AL, 36688,
United States.}
}

\date{\today}
\maketitle

\begin{abstract}
We investigate the complexities of the McKean–Vlasov optimal control problem, exploring its various formulations such as the strong and weak formulations, as well as both Markovian and non-Markovian setups within financial markets. Furthermore, we examine scenarios where the law governing the control process impacts the dynamics of options. By conceptualizing controls as probability measures on a fitting canonical space with filtrations, we unlock the potential to devise classical measurable selection methods, conditioning strategies, and concatenation arguments within this innovative framework. These tools enable us to establish the dynamic programming principle under a wide range of conditions.
\end{abstract}

\section*{Introduction:} The Black-Scholes model stands as the quintessential example of financial models, embodying key principles that have influenced subsequent developments. Its widespread adoption reflects the common adjustments and trade-offs necessary to translate theoretical concepts into practical tools for financial analysis and decision-making. Let us dive deep into the core principles of this model. The main aim of creating the Black-Scholes model is to evaluate and manage options, operating on the premise that the volatility of the underlying asset remains constant \citep{fournie1997some}. Building on this assumption, it provides an analytical formula for calculating the value of a European option, which relies on the volatility parameter of the system. However, the model's main drawback arises from the fact that in actuality, volatility fluctuates, and the market does not necessarily conform to the values predicted by this model \citep{pramanik2020optimization}.

However, practitioners in financial markets have widely embraced the Black-Scholes model, albeit with some modifications that deviate from the original theoretical expectations. These adjustments include anticipating negotiations for transactions in terms of implicit volatility rather than option prices, implementing theoretical delta hedging as recommended by the model, and recognizing that the implicit volatility tends to be higher for options further from the money and closer to maturity \citep{pramanik2021optimization,pramanik2021optimal,pramanik2023scoring}. Consequently, the Black-Scholes model retains significant practical utility, although its role differs from the primary hypotheses outlined in the original theory. While one might assume that the compromises and adaptations made in utilizing the Black-Scholes model are unique to this particular context, we argue that they offer broader insights applicable across various fields beyond finance.

So far in our initial discussion, we overlooked the inclusion of control theory. Control theory encompasses two primary forms: open loop and closed loop control. In open loop control, the decision of the policymaker remains uninfluenced by the state variable's information. Conversely, closed loop control entails the policymaker's decision being contingent upon all available state variable information \citep{pramanik2024dependence}. Integrating control variables is crucial in understanding the dynamics of financial markets. For instance, consider an entrepreneur aiming to maximize long-term profits, assuming control over advertising expenditure. If they increase spending on commercials presently, they may afford a celebrity endorsement, thereby boosting sales probability in subsequent periods \citep{pramanik2020motivation}. Conversely, heightened sales in the current period provide more funds for future commercial investment. This exemplifies closed loop control within financial markets. Furthermore, all of theses financial models are not always linear stochastic differential equations (SDEs) in practice. Therefore, the importance of non-linear SDEs come into place. One such SDE is McKean-Vlasov \citep{mckean1966class} system. The next section gives an overview of it.

\section*{Overview of McKean-Vlasov SDEs:}

McKean-Vlasov type SDEs are commonly known as nonlinear SDEs, with \emph{non-linear} highlighting the potential dependence of coefficients on the solutions' marginal distributions rather than the state variable itself. The distinctive aspect of these coefficients, even when deterministic, leads to non-local feedback effects, precluding the conventional Markov property \citep{carmona2015forward}. Incorporating the marginal distribution into the system's state could restore the Markov property but at the expense of significantly increasing the complexity by introducing a probability measure, rendering the system infinite-dimensional. While analyzing the infinitesimal generator can utilize tools developed for infinite-dimensional differential operators, conventional differential calculus, even in infinite dimensions, struggles to ensure the second component of the state process matches the statistical distribution of the first component \citep{carmona2015forward}. Nevertheless, pathwise analysis complexity remains manageable and can be effectively addressed using standard stochastic calculus tools. This points towards devising a wholly probabilistic approach to address optimal control in these nonlinear systems.

The optimal control of systems governed by McKean-Vlasov SDEs appears to be a relatively novel challenge, curiously overlooked in stochastic control literature. Typically, solving a McKean-Vlasov SDE involves employing a fixed-point argument: Initially, a set of potential marginal distribution candidates is established, and then the resulting standard SDE is solved, with the fixed-point argument entailing the requirement that the solution's marginal distributions match the initially chosen ones. Introducing stochastic control introduces an additional layer of optimization atop this fixed-point process. This formulation shares many similarities with the mean-field game (MFG) problem, as initially proposed by \cite{lasry2007mean}. In \cite{carmona2013control}, the similarities and disparities between these two problems are identified and discussed. It is notably emphasized that optimizing first and subsequently seeking a fixed point yields the solution to a mean-field game problem, while first determining the fixed point and then optimizing leads to the resolution of the optimal control for McKean-Vlasov SDEs \citep{polansky2021motif,pramanik2024estimation}. The solutions to both problems describe equilibrium states of large populations whose interactions and objectives follow mean-field patterns. The nuances between these equilibrium concepts are subtle and hinge on the optimization component's formulation within the equilibrium model. Linear quadratic models are presented in \cite{carmona2013control} to offer illustrative examples of these distinctions.

As McKean-Vlasov SDEs exhibit non-Markovian dynamics, it is logical to tackle optimization through a tailored version of the Pontryagin stochastic maximum principle, rather than forcing fit with the Hamilton-Jacobi-Bellman (HJB) framework. The stochastic Pontryagin method hinges on introducing and manipulating adjoint processes, which are solutions to adjoint backward stochastic differential equations (BSDEs). These equations involve partial derivatives of the Hamiltonian function concerning the state variable \citep{carmona2015forward}. However, in the context of McKean-Vlasov SDEs, the marginal distributions of solutions serve as full-fledged variables of the Hamiltonian function, necessitating differentiation in the pursuit of critical points. This discrepancy likely contributes to the impasse in existing literature. Until now, only dynamics contingent on moments of the marginal distributions have been explored, where differentiability with respect to the measure can be achieved through conventional calculus chain rules \citep{carmona2015forward}. The right notation of differentiability was introduced in \cite{lasry2006jeux1} and \cite{lasry2006jeux2}.

\section*{Pontryagin principle and Path integral control:}

The stochastic Pontryagin principle stands as a potent tool, yet its illuminating insights often necessitate restrictive assumptions on the models. For instance, findings in \cite{carmona2015forward} lean on a set of technical prerequisites that confine the scope of models to those characterized by coefficients primarily linear in state, control, and measure variables, alongside costs that exhibit convexity concerning the state and control variables. While these conditions may seem limiting, they are typical in applications of the stochastic Pontryagin principle to control problems. Notably, assuming the convexity of the control space is done for simplicity's sake. More extensive spaces could be accommodated, albeit at the expense of employing spike variation techniques and introducing an additional adjoint equation. However, \cite{carmona2015forward} opted against pursuing this broader generality, as it would entail added complexity in notation and technicalities, potentially obscuring the core focus of their research, especially without specific application-driven motivation \citep{pramanik2022lock,pramanik2022stochastic,pramanik2023path1,pramanik2023path}.

The essential aspect of the stochastic Pontryagin principle involves seeking within the control set for a potential minimizer of the Hamiltonian (often referred to as satisfying the Isaacs condition), while the complementary aspect suggests incorporating the minimizer's formula into both the forward equation governing dynamics and the adjoint backward equation describing adjoint processes \citep{hua2019assessing,pramanik2016tail,pramanik2021effects}. The inclusion of the minimizer in these equations establishes a significant interconnection between the forward and backward equations, effectively reducing the control problem solution to that of a forward-backward stochastic differential equation (FBSDE). Applying this approach to the current scenario entails investigating FBSDEs where the marginal distributions of solutions emerge in the coefficients of the equations. \cite{carmona2015forward} term these equations ``mean-field FBSDEs" or ``FBSDEs of McKean-Vlasov type," asserting that their study had not been undertaken prior to their research. They noted that while a general existence result was proposed in \cite{carmona2013}, one of the assumptions therein prevents its application to the linear quadratic model.

A specific variant of the stochastic Pontryagin principle, such as the Feynman-type path integral method introduced by \cite{pramanik2020optimization}, can be utilized and further developed by \cite{pramanik2024optimization} to obtain an analytical solution for the system in question \citep{pramanik2024semicooperation}. When the state variable is high-dimensional and the system dynamics are nonlinear, as seen in Merton-Garman-Hamiltonian SDEs, numerical construction of an HJB equation becomes extremely difficult. The Feynman-type path integral approach effectively overcomes this challenge of dimensionality and provides a localized analytical solution. To employ this approach, we first establish a stochastic Lagrangian for each continuous time point within the interval $s\in[0,t]$, where $t>0$. Subsequently, we divide this entire time span into $n$ equal-length subintervals and define a Riemann measure corresponding to the state variable for each subinterval. After constructing a Euclidean action function, we obtain a Schroedinger-type equation via Wick rotation. By enforcing the first-order conditions with respect to both the state and control variables, we determine the solution of the system. This methodology holds potential application in cancer research \citep{dasgupta2023frequent,hertweck2023clinicopathological,kakkat2023cardiovascular,khan2023myb,vikramdeo2023profiling,khan2024mp60}.

\section*{Further impact:}
We gain insight into the challenge of optimal control for mean-field stochastic differential equations, also known as McKean–Vlasov stochastic differential equations in academic literature. This problem presents a stochastic control scenario where the state process is driven by SDEs with coefficients influenced by the present time, the trajectories of the state process, and its distribution (or conditional distribution in cases involving shared noise). Likewise, the reward functionals can be influenced by the distribution of the state process \citep{pramanik2023optimal}.

A common limitation of the aforementioned results is their general requirement for some Markovian property of the system or its distribution, along with strong regularity assumptions on the coefficients and reward functions under consideration. This might come as a surprise to those familiar with the classical stochastic Pontryagin principle \citep{pramanik2021consensus}. In stochastic control problems, it is feasible to utilize measurable selection arguments to derive the stochastic Pontryagin principle without necessitating more than mild measurability assumptions \citep{pramanik2024bayes}. Typically, two crucial elements are needed to establish the dynamic programming principle: first, ensuring the stability of controls concerning conditioning and concatenation, and second, the measurability of the associated value function. Measurable selection arguments offer a framework to validate the conditioning, concatenation, and measurability requirements of the associated value function without relying on strong assumptions. Through the concept of relaxed control, interpreting a control as a probability measure on a canonical space, and leveraging martingale problem theory, researchers have demonstrated a stochastic Pontryagin principle using straightforward and transparent reasoning \citep{djete2022mckean}. In the future, we aim to explore scenarios where the drift and diffusion components, as well as the reward functions in option markets, can depend on the joint conditional distribution of the state process path and the closed-loop control.

\bibliographystyle{apalike}
\bibliography{bib}
\end{document}